\documentclass[12pt]{amsart}
\usepackage{amssymb,latexsym,amsmath}
\newtheorem{thm}{Theorem}

\newtheorem{lem}{Lemma}

\begin{document}
\title{On the Oppenheim's "factorisatio numerorum" function}
\author{Florian Luca, Anirban Mukhopadhyay and Kotyada Srinivas }
\address{Mathematical Institute, UNAM, Ap. Postal 61-3 (Xangari), CP 58089,
Morelia, Michoac\'an, Mexico}
\email[Florian Luca]{fluca@matmor.unam.mx}
\address{Institute of Mathematical Sciences,
CIT Campus, Taramani, Chennai 600113, India}
\email[Anirban Mukhopadhyay]{anirban@imsc.res.in}
\email[Kotyada Srinivas]{srini@imsc.res.in}

\maketitle

\section{Introduction}

Let $f(n)$ denote the number of distinct unordered factorisations of
the natural number $n$ into factors larger than $1$. For example,
$f(28)=4$ as $28$ has the following factorisations
$$
28,\ 2\cdot 14,\ 4\cdot 7,\ 2\cdot 2\cdot 7.
$$
In this paper, we address three aspects of the function $f(n)$. For
the first aspect, in \cite{CEP}, Canfield, Erd\H{o}s and Pomerance
mention without proof that the number of values of $f(n)$ that do
not exceed $x$ is $x^{o(1)}$ as $x\to\infty$. Our first theorem in
this note makes this result explicit.

For a set ${\mathcal A}$ of positive integers we put ${\mathcal
A}(x)=\{n\in {\mathcal A}: n\le x\}$.

\begin{thm}
\label{thm1} Let ${\mathcal A}=\{f(m):m\in {\mathbb N}\}$. Then
$$
\#{\mathcal A}(x)= x^{O({\log\log\log x}/{\log\log x})}, \qquad
{\text{as}}~ x\rightarrow \infty.
$$
\end{thm}

Secondly, there is a large body of literature addressing average
values of various arithmetic functions in short intervals. Our next
result gives a lower bound for the average of $f(n)$ over a short
interval.

\begin{thm}
\label{thm2} Uniformly for $x\ge 1$ and $y>e^{e^{e}}$, we have
$$
\frac{1}{y} \sum_{x\leq n\leq x+y} f(n) \geq \exp
\left(\left(\frac{4}{\sqrt{2e}} +  O\left(\frac{(\log\log\log
y)^2}{\log\log y}\right)\right)  \frac{{\sqrt{\log y}}}{\log\log y}
\right).
$$
\end{thm}

Finally, there are also several results addressing the behavior of
positive integers $n$ which are multiples of some other arithmetic
function of $n$. See, for example, \cite{CK}, \cite{EP}, \cite{S1}
and \cite{S2} for problems related to counting positive integers $n$
which are divisible by either $\omega(n),~\Omega(n)$ or $\tau(n)$,
where these functions are the number of distinct prime factors of
$n$, the number of total prime factors of $n$, and number of
divisors of $n$, respectively. Our next and last result gives an
upper bound on the counting function of the set of positive integers
$n$ which are multiples of $f(n)$.

\begin{thm}
\label{thm3} Let $ \mathcal{B} = \left\{ n : f(n) \mid n \right\}$.
Then
$$
\# \mathcal{B}(x)= \frac{x}{(\log x)^{1 + o(1)}} \ ,\qquad
{\text{as}}~ x\rightarrow \infty.
$$
\end{thm}

\section{Preliminaries and lemmas}
\label{sec:Intro}

The function $f(n)$ is related to various partition functions. For
example, $f(2^n)=p(n)$, where $p(n)$ is the number of partitions of
$n$. Furthermore, $f(p_1p_2\cdots p_k)=B_k$, where $B_k$ is the
$k$th Bell number which counts the number of partitions of a set
with $k$ elements in nonempty disjoint subsets. In general,
$f(p_1^{\alpha_1}p_2^{\alpha_2}\cdots p_k^{\alpha_k})$ is the number
of partitions of a multiset consisting of $\alpha_i$ copies of
$\{i\}$ for each $i=1,\ldots,k$. Throughout the paper, we put $\log
x$ for the natural logarithm of $x$. We use $p$ and $q$ for prime
numbers and $O$ and $o$ for the Landau symbols.

The following asymptotic formula for the $k$th Bell number is due to
de Brujin \cite{DBR}.

\begin{lem}
\label{bell}
$$
\frac{\log B_k}{k}=\log k-\log\log k-1
  +\frac{\log\log k}{\log k}+\frac{1}{\log k}+O\left(\frac{(\log\log k)^2}{(\log k)^2}\right).
  $$
\end{lem}
We also need the Stirling numbers of the second kind $S(k,l)$ which
count the number of partitions of a $k$ element set into $l$
nonempty disjoint subsets. Clearly,
\begin{equation}
\label{eq:Bk}
B_k=\sum_{l=1}^k S(k,l).
\end{equation}
We now formulate and prove a few lemmas about the function $f(n)$
which will come in handy later on.

The first lemma is an easy observation, so we state it without
proof.
\begin{lem}
\label{div} If $ a \mid b $, then $f(a) \leq f(b)$.
\end{lem}

We let $p_n$ denote the $n$th prime number and $\alpha_1(n)$ denote
the maximal exponent of a prime appearing in the prime factorization
of $n$. Let $n$ be a positive integer with prime factorization
$$
n=q_1^{\alpha_1}q_2^{\alpha_2}\cdots q_k^{\alpha_k},
$$
where $q_1,\ldots,q_k$ are distinct primes and
$\alpha_1(n):=\alpha_1\ge \alpha_2\ge\cdots \ge \alpha_k$. We put
$n_0=p_1^{\alpha_1}p_2^{\alpha_2}\cdots p_k^{\alpha_k}$, and observe
that $f(n)=f(n_0)$. This observation will play a crucial role in the
proof of Theorem \ref{thm1}.

The following lemma gives upper bounds for $\alpha_1(n)$ and
$\omega(n)$ when $f(n)\le x$.
\begin{lem}
\label{alphat} Let $n = q_1^{{\alpha}_1}q_2^{{\alpha}_2}\cdots
{q_k}^{{\alpha}_k}$, where $\alpha_1 \ge \alpha_2 \ge \cdots \ge
\alpha_k$ and $f(n) \leq x$. Then
\begin{enumerate}
\item[(i)] $ {\alpha}_1 = O \big( (\log x)^2 \big);$
\item[(ii)] $ k = \omega(n)= O \big({ \log x}/{\log\log x} \big).$
\end{enumerate}
\end{lem}

\proof It follows from Lemma \ref{div} that
$$f(n)\ge f(q_1^{{\alpha}_1})=p(\alpha_1).$$
Using the following asymptotic formula for $p(n)$
due to Hardy and Ramanujan \cite{HR}
\begin{equation}\label{pn}
p(n)\sim \frac{1}{4n\sqrt 3}\exp(\pi\sqrt{2n/3}),
\end{equation}
we conclude that $\exp(c\sqrt{\alpha_1})\le x$ holds with some
constant $c>0$. Hence, (i) follows. In order to prove (ii), let
again $n_0=p_1p_2\cdots p_k$. By Lemma \ref{div}, we have $f(n_0)\le
f(n)\le x$. Furthermore, $f(n_0)=B_k$. It now follows from Lemma
\ref{bell}, that
$$
\exp\big((1+o(1))k\log k)=B_k\le x,
$$
as $k\to\infty$, yielding
$$
k= O \left( \frac{ \log x}{\log\log x} \right),
$$
which completes the proof of the lemma.
\qed

\medskip

Recall that the M\"obius function $\mu(m)$ of the positive integer
$m$ is $(-1)^{\omega(m)}$ if $m$ is squarefree and $0$ otherwise.

For a positive integer $k$ and positive real numbers $A\le B$ we let
$$
\mathcal{M}_{k,A,B} = \left\{m: \mu(m)\ne 0,~ \omega
(m)=k,~{\text{\rm if}}~p\mid m~{\text{\rm then}}~ p\in [A,B]
\right\}.
$$
We also put
$$
S_{A,B}=\sum_{A\le p\le B}\frac{1}{p}.
$$

\begin{lem}
\label{reciprocal} Uniformly in $A\ge 2,~B\ge 3$ and $k\ge 2$, we
have
$$
\sum_{m\in \mathcal{M}_{k,A,B}}\frac{1}{m} \ge
\left(1+O\left(\frac{k^2}{S_{A,B}^2A\log A}\right)\right)
\frac{1}{k!}S_{A,B}^k.
$$
\end{lem}

\proof We omit the dependence of the subscripts in order to simplify
the presentation. It is not hard to see that
\begin{equation}\label{lemma4main}
\sum_{m\in \mathcal{M}}\frac{1}{m} \geq \frac{1}{k!}\Big(\sum_{A\leq
p \leq B}\frac{1}{p} \Big)^k - \sum_{A\le p\le
B}\frac{1}{p^2}\frac{1}{(k-2)!} \Big(\sum_{A\leq p \leq
B}\frac{1}{p}\Big)^{k-2}.
\end{equation}
Indeed, if $ m = {q_1}^{\alpha_1}\cdots {q_s}^{\alpha_s}$, with
$\alpha_1 \geq 2$ and $\alpha_1+\cdots +\alpha_s =k$, then, by
unique factorization, in the first sum on the right hand side of
inequality (\ref{lemma4main}), the number $1/m$ appears with
coefficient
$$
\frac{1}{k!} \left( \frac{k!}{\alpha_1!\cdots \alpha_s! }\right) =
\frac{1}{{\alpha_1}! \cdots {\alpha_s}!},
$$
while in the second sum in the right hand side of the inequality
(\ref{lemma4main}), the number $1/m$ appears with coefficient
\begin{equation*}
\sum_{\substack{1\le i\le s\\ \alpha_i\ge 2}}\frac{1}{(k-2)!} \left(
\frac{(k-2)!}{\alpha_1!\cdots (\alpha_i-2)!\cdots \alpha_s!}\right)
> \frac{1}{\alpha_1 ! \cdots \alpha_s !}.
\end{equation*}
This establishes inequality (\ref{lemma4main}). Using this
inequality, we get
\begin{eqnarray*}
\sum_{m \in \mathcal{M}}\frac{1}{m} \geq \frac{S^k}{k!}
&-&    \frac{1}{(k-2)!}S^{k-2} \sum_{A\le p\le B} \frac{1}{p^2}\\
&\ge & \frac{S^k}{k!}
\left( 1 - \frac{k^2}{S^2} \sum_{p\geq A} \frac{1}{p^2}\right).
\end{eqnarray*}
An argument involving the Prime Number Theorem and partial summation
gives
$$
\sum_{p\geq A} \frac{1}{p^2}=O\left(\frac{1}{A\log A}\right).
$$
Hence,
$$ \sum_{m \in \mathcal{M}}\frac{1}{m} \geq \frac{S^k}{k!}
\left(1+O\left(\frac{k^2}{S^2 A\log A}\right)\right).
$$
This completes the  proof of the lemma.

\section{Proofs of the theorems}

\subsection{Proof of Theorem \ref{thm1}}

For a positive integer $n$, we let again $n_0$ and $\alpha_1(n)$ be
the functions defined earlier. We let ${\mathcal
A}(x)=\{m_1,\ldots,m_t\}$ be such that $m_1<m_2<\cdots <m_t$ and let
${\mathcal N}=\{n_1,\ldots,n_t\}$ be positive integers such that
$n_i$ is minimal among all positive integers $n$ with $f(n)=m_i$ for
all $i=1,\ldots,t$. It is clear that if $n\in {\mathcal N}$, then
$n=n_0$. Since $\#{\mathcal A}(x)=t=\#{\mathcal N}$, it suffices to
bound the cardinality of ${\mathcal N}$.

We partition this set as $\mathcal N=\mathcal N_1\cup\mathcal
N_2\cup\mathcal N_3$, where
$$
\mathcal N_1=\left\{n\in\mathcal N: \alpha_1(n)\le \log\log x
\right\},
$$
$$
\mathcal N_2=\left\{n\in\mathcal N: \omega(n) \le \frac{\log
x}{(\log\log x)^2}\right\},
$$
and
$$
\mathcal N_3=\mathcal N \backslash \mathcal N_1\cup\mathcal N_2.
$$

If $n\in \mathcal N_1$, then $n$ has at most $O(\log x/\log\log x)$
prime factors (by Lemma \ref{alphat}), each one appearing at an
exponent less than $\log\log x$.

Therefore,
\begin{eqnarray}
\label{eq:N1}
\#\mathcal N_1 & = & (\log\log x)^{O({\log
x}/{\log\log x})} =\exp \left(O \left(\frac{\log x \log\log\log
x}{\log\log x}\right)
\right)\nonumber\\
& = & x^{O\left(\frac{\log\log\log x}{\log\log x}\right)}
\end{eqnarray}
as $x\to\infty$.

Next, we observe that an integer in $\mathcal N_2$ has at most $\log
x/(\log\log x)^2$ prime factors, each appearing at an exponent
$O((\log x)^2)$ (by Lemma \ref{alphat}). Thus,
\begin{eqnarray}
\label{eq:N2} \# \mathcal{N}_2 & \leq & \left(O\left(\log x
\right)^2 \right)^{{\log x}/{(\log\log x)^2}} = \exp \left(
\frac{\left(2+o(1)\right)\log x}{\log\log x}\right)\nonumber\\
&  = & x^{o\left(\frac{\log\log\log x}{\log\log x}\right)}
\end{eqnarray}
as $x\to\infty$.

Finally, let $n\in \mathcal{N}_3$, and write it as
$$
n=p_1^{{\alpha}_1}\cdots
{p_i}^{\alpha_i}{p_{i+1}}^{\alpha_{i+1}}\cdots{p_k}^{\alpha_k},
$$
where we put
$$
i: = \max \{ j \leq k : \alpha_j \geq y \},
$$
where $y=\log\log x/\log\log\log x$.

Observe that the divisor ${p_{i+1}}^{\alpha_{i+1}}\cdots
{p_t}^{\alpha_t} $ of $n$ can be chosen in at most
\begin{equation}
\label{eq:ytok} (y+1)^{k}= (y+1)^{O({\log x}/{\log\log x})} =
 \exp \left(O\left(\frac{\log x \log \log\log x}{\log\log x} \right)\right)
\end{equation}
ways. Furthermore, by Lemma \ref{alphat}, we trivially have that
$n'=p_1^{\alpha_1}\cdots p_i^{\alpha_i}$ can be chosen in at most
$$
\left(O\left((\log x)^2\right)\right)^i=\exp\left(O(i\log\log
x)\right).
$$
Thus, putting ${\mathcal N}_4$ for the subset of ${\mathcal N}_3$
such that $i\le \log x/(\log\log x)^2$, we get that
\begin{equation}
\label{eq:N4} \#{\mathcal N}_4\le \exp\left(O\left(\frac{\log
x}{\log\log x}\right)\right).
\end{equation}
From now on, we look at $n\in {\mathcal N}_5={\mathcal
N}_3\backslash {\mathcal N}_4$.

For each $t$, we let $k_t$ be such that $S(t,k_t)$ is maximal among
the numbers $S(t,k)$ for $ k= 1,\ldots,t$. By formula \eqref{eq:Bk},
the definition of $k_t$, and Lemma \ref{bell}, we have that
$$
S(t,k_t) \geq \frac{B_t}{t}=\frac{\exp((1+o(1))t\log
t)}{t}=\exp((1+o(1))t\log t)
$$
as $t\to\infty$. We now claim that
$$
f(n) \geq f(n') \ge f((p_1\cdots p_i)^y) \ge
\frac{S(i,k_i)^y}{(yk_i)!}.
$$
The first three inequalities follow immediately from Lemma
\ref{div}, so let us prove the last one.

Note that $S(i,k_i)$ counts the number of factorizations of
$p_1p_2\cdots p_i$ in precisely $k_i$ factors. Therefore,
$\big(S(i,k_i)\big)^y$ counts the number of factorizations of $
(p_1p_2\cdots p_i)^y$ into $k_iy$ square-free factors, where we
count each such factorization at most $(k_iy)!$ times. This
establishes the claim.

Since $i$ tends to infinity for $n\in {\mathcal N}_5$, we get that
$$
S(i,k_i)^y \geq \exp \left( (1+o(1))yi\log i\right).
$$
Furthermore, we trivially have
$$
(k_iy )! \leq (k_iy)^{k_iy}=\exp\left (k_iy\log (k_iy)\right).
$$
Thus,
\begin{equation}
\label{eq:Stir} f(n)\ge \frac{S(i,k_i)^y}{(k_iy)!}\ge
\exp\left((1+o(1))yi\log i-k_i y\log(k_iy)\right)
\end{equation}
as $x\to\infty$. We next show that for our choices of $y$ and $i$ we
have
$$
k_i y\log (k_iy)=o(yi\log i)\qquad {\text{\rm as}}\qquad x\to\infty.
$$
Indeed, using the fact
$$k_i=(1+o(1)){\displaystyle{\frac{i}{\log i}}}\qquad {\text{\rm
as}}\qquad  i\to\infty
$$
(see, for example, \cite{CP}), we see that the above condition is
equivalent to
$$
\log y= o((\log i)^2),
$$
which holds as $x\to\infty$ because $y=\log\log x/\log\log\log x$
and $i> \log x/(\log\log x)^2$. Now the inequality $f(n)\le x$
together with \eqref{eq:Stir} and the fact that $\log i\ge
(1+o(1))\log \log x$ implies that
\begin{equation}
\label{eq:i} i\le (1+o(1))\frac{\log x}{y\log\log x}\qquad
{\text{\rm as}}\quad x\to\infty,
\end{equation}
therefore $n'$ can be chosen in at most
\begin{equation}\label{n'}
\left(O\left((\log x)^2\right)\right)^i\le \left(O\left((\log
x)^2\right)\right)^{(1+o(1))\frac{\log x}{y\log\log x}}
  =\exp\left((2+o(1))\frac{\log x}{y}\right)
\end{equation}
ways. As we have already seen at \eqref{eq:ytok}, the complementary
divisor $n/n'=p_{i+1}^{\alpha_{i+1}}\cdots p_t^{\alpha_t}$ of $n$
can be chosen in at most
\begin{equation}\label{n/n'}
\exp\left(O\left(\frac{\log x\log\log\log x}{\log\log
x}\right)\right)
\end{equation}
ways. Thus, the total number of choices for $n$  in ${\mathcal N}_5$
is
\begin{eqnarray}
\label{eq:n3} \#{\mathcal N}_5 & \le & \exp\left(O\left(\frac{\log
x}{y}+\frac{\log x\log y}{\log\log
x}\right)\right)\nonumber\\
& = & \exp\left(O\left(\frac{\log x\log\log\log x}{\log\log
x}\right)\right).
\end{eqnarray}
Hence, from estimates \eqref{eq:N4} and \eqref{eq:n3} we get
\begin{equation}
\label{eq:N3} \#{\mathcal N}_3\le \#{\mathcal N}_4+\#{\mathcal
N}_5\le x^{O(\log\log\log x/\log\log x)}.
\end{equation}
From estimates \eqref{eq:N1}, \eqref{eq:N2} and \eqref{eq:N3},  we
finally get
$$
\#{\mathcal N}\le \#{\mathcal N}_1+\#{\mathcal N}_2+\#\mathcal N_3
\le x^{O(\log\log\log x/\log\log x)},
$$
which completes the
proof of the theorem.

\subsection{Proof of Theorem \ref{thm2}}
We assume that $y$ is as large as we wish otherwise there is nothing
to prove. Let $s=\lfloor 3\log\log y\rfloor$. Let
$$
\mathcal{N} = \left\{n \in(x,x+y): n~{\text{\rm has}}~k+j~{\text{\rm
prime~factors~in}}~[A,B],~0\le j\le s-1\right\},
$$
with the parameters $A=k^2, B=y^{1/{(k+s+1)}}$, where we take $k\in
[c_1\sqrt{\log y},~c_2{\sqrt{\log y}}]$, and $0<c_1<c_2$ are two
constants to be made more precise later. We will spend some time
getting a lower bound on the cardinality of ${\mathcal N}$. For
this, observe that for each $n\in {\mathcal N}$ there is a
squarefree number $m$ with exactly $k$ distinct prime factors in
$[A,B]$ such that $m\mid n$. Clearly, $m\le y^{k/(k+s+1)}$. Fix such
an $m$ and put ${\mathcal N}_m$ for the set of multiples of $m$ in
${\mathcal N}$. To get a lower bound on $\#{\mathcal N}_m$, observe
first that the number of multiples of $m$ in $(x,x+y)$ is
$$
\ge \left\lfloor\frac{y}{m}\right\rfloor\ge
\frac{y}{m}-1=\frac{y}{m}\left(1+O\left(\frac{m}{y}\right)\right)=\frac{y}{m}\left(1+
O\left(\frac{1}{\log y}\right)\right).
$$
Of course, not all such numbers are in ${\mathcal N}_m$ since some
of them might have more than $k+s-1$ distinct prime factors in
$[A,B]$. We next get an upper bound for the number of such ``bad"
multiples $n$ of $m$. For each such bad $n$, there exists a number
$m_1$ having $s$ prime factors in $[A,B]$ and coprime to $m$ such
that $mm_1\mid n$. Note that $mm_1\le y^{(k+s)/(k+s+1)}<y$. For
fixed $m$ and $m_1$, the number of such positive integers in
$(x,x+y)$ is
$$
\le \left\lfloor\frac{y}{mm_1}\right\rfloor+1\le \frac{2y}{mm_1}.
$$
Summing up over all possibilities for $m_1$, we get that the number
of such $n$ is
$$
\le \frac{2y}{m}\sum_{\substack{m_1\in {\mathcal M}_{s,A,B}\\
(m_1,m)=1}}\frac{1}{m_1}\le \frac{2y}{ms!}\left(\sum_{A\le p\le
B}\frac{1}{p}\right)^s=\frac{2yS^s}{ms!},
$$
where we put
$$
S:=\sum_{A\le p\le B}\frac{1}{p}.
$$
Observe that, by Mertens's formula, we have
\begin{eqnarray*}
S & = & \left(\log\log B+c_0\right)-\left(\log\log
A+c_0\right)+O\left(\frac{1}{\log
A}\right)\\
& = & \log\log y-\log(k+s+1)-\log\log k-\log 2+O\left(\frac{1}{\log
k}\right)\\
& = & \log\log y-\log k-\log\log k-\log 2+O\left(\frac{1}{\log
k}+\frac{s}{k}\right).
\end{eqnarray*}
As far as errors go, note that since $s=3\log\log y+O(1)$ and
$k\asymp {\sqrt{\log y}}$, we have that
$$
\frac{s}{k}\ll \frac{\log k}{k}\ll \frac{1}{\log k}.
$$
Furthermore, $S=(1/2+o(1))\log\log y$ as $y\to\infty$, therefore for
$y>y_0$ we have that $S<s/3$. We record that
\begin{equation}
\label{eq:Mert} S=\log\log y-\log k-\log\log k-\log
2+O\left(\frac{1}{\log k}\right).
\end{equation}
The above arguments show that
$$
\#{\mathcal N}_m \ge
\frac{y}{m}\left(1-\frac{2S^s}{s!}+O\left(\frac{1}{\log
y}\right)\right).
$$
From the elementary estimate $s!\ge (s/e)^s$, we get
$$
\frac{2S^s}{s!}\ll \left(\frac{Se}{s}\right)^s\ll
\left(\frac{e}{3}\right)^s
$$
and the last number above is $<1/3$ if $y$ is sufficiently large.
Hence, the inequality
$$
\#{\mathcal N}_m\ge \frac{y}{2m}
$$
holds uniformly in squarefree integers $m$ having $k$ distinct prime
factors all in $[A,B]$. We now sum over $m$ and use Lemma
\ref{reciprocal} to get that
\begin{equation}
\label{eq:Nm}
\sum_{m\in {\mathcal M}_{k,A,B}}\#{\mathcal N}_m\ge
\frac{y}{2}\sum_{m\in {\mathcal M}_{k,A,B}}\frac{1}{m}\gg
\frac{yS^k}{k!}\left(1+O\left(\frac{k^2}{S^2A\log
A}\right)\right)\gg \frac{yS^k}{k!}
\end{equation}
for large $y$, because $A=k^2$, therefore the expression
$k^2/(S^2A\log A)$ is arbitrarily small if $y$ is large. Next let us
note that if $n\in {\mathcal N}$, then $n$ has $k+j$ distinct prime
factors in $[A,B]$ for some $j=0,1,\ldots,s-1$. Thus, the number of
possibilities for $m\mid n$ in ${\mathcal M}_{k,A,B}$ is
$$
\binom{k+j}{k}\le
\binom{k+s}{s}<\left(e+\frac{ek}{s}\right)^s=\exp(O((\log\log
y)^2)).
$$
Here, we used again the fact that $s!\ge (s/e)^s$. In particular,
the sum on the left of \eqref{eq:Nm} counts numbers $n\in {\mathcal
N}$ and each number is counted at most $\exp(O((\log\log y)^2))$
times. Hence, dividing by this number we get a lower bound on
$\#{\mathcal N}$ which is
$$
\#{\mathcal N}\ge \frac{yS^k}{k!}\exp(O((\log\log y)^2)).
$$
If $n\in \mathcal N$, then there is an $m\in\mathcal M$ such that
$m\mid n$. It now follows, from Lemma \ref{div}, that $f(n) \geq
f(m) \geq B_k.$ Thus,
$$\frac{1}{y}\sum_{x\leq n \leq x+y }f(n)
\ge \frac{1}{y}\sum_{n\in\mathcal N} f(n) \ge
\frac{1}{y}B_k\#\mathcal N \ge  \frac{B_kS^k}{k!}\exp(O((\log\log
y)^2)).
$$
We now maximize $B_k{S^k}/{k!}$ by choosing $k$ appropriately versus
$y$. Using Stirling's formula
$$
k! \sim \left( \frac{k}{e} \right)^k (2 \pi k)^{1/2}
$$
to estimate $k!$, Lemma \ref{bell} as well as estimate
\eqref{eq:Mert}, we get
\begin{equation}
\label{eq:function} \frac{B_kS^k}{k!}\exp(O((\log\log y)^2)) =
\exp\left(h(k)+O\left(\frac{k(\log\log k)^2}{(\log
k)^2}\right)\right),
\end{equation}
where the function $h(k)$ is
\begin{eqnarray*}
h(k) &=& k\log (\log\log y-\log k-\log\log k-\log 2)\\
& - & k\log\log k
     + k\frac{\log\log k}{\log k}+\frac{k}{\log k}.
\end{eqnarray*}
The error term under the exponential in formula \eqref{eq:function}
comes from the estimate given by Lemma \ref{bell} on $B_k$, estimate
\eqref{eq:Mert} which tells us that
\begin{eqnarray*}
k\log S & = & k\log\left(\log\log y-\log k-\log\log k-\log
2+O\left(\frac{1}{\log k}\right)\right)\\
& = & k\log(\log\log y-\log k-\log\log k-\log
2)+O\left(\frac{k}{(\log k)^2}\right),
\end{eqnarray*}
because $\log\log y-\log k-\log\log k-\log 2\asymp \log k$ for our
choice of $k$ versus $y$, as well as the fact that $(\log\log
y)^2\ll k(\log\log k)^2/(\log k)^2$, again by our choice of $k$
versus $y$.

We now choose
$$
k=\left\lfloor \frac{1}{\sqrt{2e}} (\log y)^{1/2}\right\rfloor.
$$
Note that with $c_1=1/4$ and $c_2=1/2$ we indeed have that $k\in
[c_1(\log y)^{1/2},~c_2(\log y)^{1/2}]$, as promised. Then,
\begin{eqnarray*}
k & = & \frac{1}{\sqrt{2e}} (\log y)^{1/2}+O(1);\\
\log k & = & \frac{1}{2}\log\log
y-\log({\sqrt{2e}})+O\left(\frac{1}{{\sqrt{\log y}}}\right);\\
\frac{1}{\log k} & = & \frac{2}{\log\log
y}+O\left(\frac{1}{(\log\log y)^2}\right).
\end{eqnarray*}
In particular,
\begin{eqnarray*}
\log\log y & - & \log k  -  \log\log k-\log 2\\
& = &  \frac{1}{2}\log\log y+\log({\sqrt{2e}}/2)-\log\log
k+O\left(\frac{1}{\sqrt{\log
y}}\right)\\
& = & \left(\frac{1}{2}\log\log y-\log({\sqrt{2e}})\right)-\log\log
k+1+O\left(\frac{1}{\sqrt{\log y}}\right)\\
& = & \log k-(\log\log k-1)+O\left(\frac{1}{\sqrt{\log y}}\right),
\end{eqnarray*}
so that
\begin{eqnarray*}
\log(\log\log y & - & \log k-\log\log k-\log 2)\\
& = & \log\left(\log k-(\log\log k-1)+O\left(\frac{1}{\sqrt{\log
y}}\right)\right) \\
& = & \log(\log k-(\log\log
k-1))+O \left(\frac{1}{k{\sqrt{\log y}}}\right)\\
& = & \log(\log k-(\log\log k-1))+O\left(\frac{1}{\log y}\right).
\end{eqnarray*}
Thus,
\begin{eqnarray*}
k\log (\log\log y & - & \log k-\log\log k-\log 2)-  k\log\log k\\
& =
&
 k\log\left(\frac{\log k-(\log\log
k-1)}{\log k}\right)\left(1+O\left(\frac{1}{{\log y}}\right)\right)\\
& = & k\log\left(1-\frac{\log\log k-1}{\log
k}\right)+O\left(\frac{1}{\sqrt{\log y}}\right)\\
& = & -\frac{k(\log\log k-1)}{\log k}+O\left(\frac{k(\log\log
k)^2}{(\log k)^2}+\frac{1}{k}\right)\\
& = & -\frac{k\log\log k}{\log k}+\frac{k}{\log
k}+O\left(\frac{k(\log\log k)^2}{(\log k)^2}\right).
\end{eqnarray*}
It now follows immediately that
\begin{eqnarray*}
h(k) & = & k\log(\log\log y  -  \log k-\log\log k-\log 2)-
k\log\log k \\
 & + & \frac{k\log\log k}{\log k}+\frac{k}{\log k}
  = \frac{2k}{\log k}+O\left(\frac{k(\log\log k)^2}{(\log
 k)^2}\right).
 \end{eqnarray*}
One can in fact check that the above estimate is the maximum of
$h(k)$ as a function of $k$ when $y$ is fixed. We will not drag the
reader through this computation. Comparing the above estimate with
\eqref{eq:function}, we get that
\begin{eqnarray*}
&& \frac{B_kS^k}{k!}\exp(O((\log\log y)^2)) \ge
\exp\left(\frac{2k}{\log k}+O\left(\frac{k(\log\log k)^2}{(\log
k)^2}\right)\right)\\
&& =  \exp\left(\frac{4}{\sqrt{2e}}\frac{(\log y)^{1/2}}{\log\log
y}\left(1+O\left(\frac{(\log\log\log y)^2}{\log\log
y}\right)\right)\right).
\end{eqnarray*}
We thus get that
\begin{eqnarray*}
\frac{1}{y}\sum_{x\leq n \leq x+y }f(n)  & \ge &
\frac{B_kS^k}{k!}\exp(O((\log\log y)^2))\\
& \ge & \exp \left(\left(\frac{4}{\sqrt{2e}} +
  O\left(\frac{(\log\log\log y)^2}{\log\log y}\right)\right)
  \frac{{\sqrt{\log y}}}{\log\log y} \right),
\end{eqnarray*}
which is what we wanted.

\subsection{Proof of Theorem \ref{thm3}}

We observe that primes are in $\mathcal A$ as $f(p)=1$ for all prime
$p$. Thus,
$$
\#\mathcal A (x) \gg \frac{x}{\log x}.
$$
This completes the lower bound part of the theorem. To obtain the
upper bound, we cover the set $\mathcal A(x)$ by three subsets
${\mathcal A}_1,~{\mathcal A}_2$ and ${\mathcal A}_3$ as follows:
$$
{\mathcal A}_1=\left\{n\le x\ : \ \Omega(n)>10\log\log x\right\},
$$
$$
{\mathcal A}_2=\left\{n\le x\ : \ \omega(n)<\frac{\log\log
x}{\log\log\log x}\right\},
$$
and
$$
{\mathcal A}_3 = \left\{ n \leq x : n \equiv 0 \pmod{f(n)},
~n\not\in {\mathcal A}_1\cup {\mathcal A}_2\right\}.$$

We recall the following bound
$$
\# \left\{ n \leq x : \Omega (n) = k \right\} \ll \frac{kx}{2^k}
$$
valid uniformly in $k$ (see, for example, Lemma 13 in \cite{LP}).
Using the above estimate, we get
\begin{equation}
\label{eq:A1}
\#\mathcal{A}_1 \leq x \sum_{k > 10 \log\log x}
\frac{k}{2^k} \ll\frac{ x \log\log x}{2^{10\log\log
x}}=o\left(\frac{x}{\log x}\right)
\end{equation}
as $x\to\infty$. To find an upper bound for $\mathcal A_2$, we use
the Hardy-Ramanujan bounds (see \cite{HR})
$$
\#\left\{n \leq x : \omega (n) = k \right\} \ll \frac{x
\left(\log\log x+c_1\right)^{k-1}}{\log x (k-1)!}
$$
with some positive constant $c_1$. Using the elementary estimate
$m!\ge (m/e)^m$ with $m=k-1$, we get
$$
\#\left\{n \leq x : \omega (n)=k\right\}
  \ll \frac{x}{\log x} \left(\frac{e\log\log
  x+c_2}{k-1}\right)^{k-1},
$$
where $c_2=ec_1$. The right hand side is an increasing function of
$k$ in our range for $k$ versus $x$ when $x$ is large. Since
$k<(\log\log x)/(\log\log\log x)$, we deduce that
\begin{equation}
\label{eq:A2}
\# \mathcal{A}_2 \ll  \frac{x}{\log x}
\left(O(\log\log\log x)\right)^{{\log\log x}/{\log\log\log
x}}=\frac{x}{(\log x)^{1+o(1)}}
\end{equation}
as $x\to\infty$.

Finally, we estimate $\mathcal A_3$. Each $n\in\mathcal A_3$ can be
written as
$$
n=q_1^{\alpha_1}q_2^{\alpha_2}\cdots q_k^{\alpha_k},
$$
where $q_1,\cdots,q_k$ are distinct primes, $\alpha_1\ge \alpha_2\ge
\cdots \ge \alpha_k$, $\alpha_1+\alpha_2+\cdots+\alpha_k\le
10\log\log x$ and $k>K:=\lfloor \log\log x/\log\log\log x\rfloor$.
Let ${\mathcal T}$ be the set of all such tuples
$(k,\alpha_1,\ldots,\alpha_k)$. For each such $n$, we have that
\begin{eqnarray*}
f(n) & \ge & B_K \ge \exp((1+o(1))K \log K)
\ge \exp((1+o(1)) \log\log x)\\
& = & (\log x)^{1+o(1)}.
\end{eqnarray*}
The number of tuples $(k,\alpha_1,\ldots,\alpha_k)$ satisfying the
above conditions is at most
$$
\#{\mathcal T}\ll \log\log x\sum_{n\le 10\log\log x}p(n),
$$
where again $p(n)$ is the partition function of $n$. Using estimate
(\ref{pn}), we get that the cardinality of ${\mathcal T}$ is at most
$$
\#{\mathcal T}\ll (\log\log x)^2\exp(O(\sqrt{\log\log x}))=(\log
x)^{o(1)}\qquad {\text{\rm as}}\quad x\to\infty.$$ Thus,
\begin{equation}
\label{eq:A3}
\# \mathcal A_3  \ll
\sum_{(k,\alpha_1,\ldots,\alpha_k)\in {\mathcal
T}}\frac{x}{f({p_1}^{\alpha_1}\cdots {p_k}^{\alpha_k})}\ll
\frac{x\#{\mathcal T}}{B_K}=\frac{x}{(\log x)^{1+o(1)}}
\end{equation}
as $x\to\infty$. Now inequalities \eqref{eq:A1}, \eqref{eq:A2} and
\eqref{eq:A3} yield the desired upper bound and complete the proof.

\end{document}